%% file: b4_1.tex
\documentstyle[12pt, epsf, epsfig, amssymb, amstext, amstex,amsthm, amsfonts]{article}
\input xy
\input amssymb.sty
\xyoption{all}
\begin{document}

\newtheorem{thm}{Theorem}[section]
\newtheorem{prop}[thm]{Proposition}
\newtheorem{cor}[thm]{Corollary}
\newtheorem{lem}[thm]{Lemma}
\newtheorem{conj}[thm]{Conjecture}
\newtheorem{exa}[thm]{Example}
\newtheorem{defn}[thm]{Definition}
\newtheorem{clm}[thm]{Claim}
\newtheorem{eex}[thm]{Exercise}
\newtheorem{obs}[thm]{Observation}
\newtheorem{note}[thm]{Notation}
 
\newcommand{\ben}{\begin{enumerate}}
\newcommand{\een}{\end{enumerate}}
\newcommand{\blem}{\begin{lem}}
\newcommand{\elem}{\end{lem}}
\newcommand{\bcl}{\begin{clm}}
\newcommand{\ecl}{\end{clm}}
\newcommand{\bthm}{\begin{thm}}
\newcommand{\ethm}{\end{thm}}
\newcommand{\bpr}{\begin{prop}}
\newcommand{\epr}{\end{prop}}
\newcommand{\bco}{\begin{cor}}
\newcommand{\eco}{\end{cor}}
\newcommand{\bcon}{\begin{conj}}
\newcommand{\econ}{\end{conj}}
\newcommand{\bde}{\begin{defn}}
\newcommand{\ede}{\end{defn}}
\newcommand{\bex}{\begin{exa}}
\newcommand{\eexa}{\end{exa}}
\newcommand{\bexe}{\begin{exe}}
\newcommand{\eexe}{\end{exe}}
\newcommand{\bobs}{\begin{obs}}
\newcommand{\eobs}{\end{obs}}
\newcommand{\bnote}{\begin{note}}
\newcommand{\enote}{\end{note}}

\newcommand{\fg}{\Pi _1(D-K,u)}
\newcommand{\Z}{{\Bbb Z}}
\newcommand{\C}{{\Bbb C}}
\newcommand{\R}{{\Bbb R}}
\newcommand{\Q}{{\Bbb Q}}
\newcommand{\F}{{\Bbb F}}
\newcommand{\N}{{\Bbb N}}

\newcommand{\fnref}[1]{~(\ref{#1})}
\newenvironment{emphit}{\begin{itemize} \em}{\end{itemize}}
\begin{center}
\Large{\bf {Properties of Hurwitz Equivalence in the Braid Group of Order $n$}}\\ 
\vspace{7mm}
\large{T. Ben-Itzhak and M. Teicher}
\footnote{Partially supported by the Emmy Noether Research Institute for Mathematics
and the Minerva Foundation of Germany and to the Excellency Center "Group
Theoretic Methods in the Study of Algebraic Varieties"  of the Israel Science Foundation.\\}

\end{center}


\vspace{15mm}
ABSTRACT. In this paper we prove certain Hurwitz equivalence properties in $B_n$. Our main result is that every two Artin's factorizations of $\Delta _n ^2$ of the form $H_{i_1} \cdots H_{i_{n(n-1)}}, \quad F_{j_1} \cdots F_{j_{n(n-1)}}$ (with $i_k , j_k \in \{ 1,...,n-1 \}$), where $\{ H_1,...,H_{n-1} \}, \{ F_1,...,F_{n-1} \}$ are frames, are Hurwitz equivalent. This theorem is a generalization of the theorem we have proved in \cite{Tz}, using an algebraic approach unlike the proof in \cite{Tz} which is geometric.\\
The application of the result is applied to the classification of algebraic surfaces up to deformation. It is already known that there exist surfaces that are diffeomorphic to each other but are not deformations of each other (Manetti example). We construct a new invariant based on a Hurwitz equivalence class of factorization, to distinguish among diffeomorphic surfaces which are not deformations of each other. The precise definition of the new invariant can be found in \cite{KuTe} or \cite{Te}.
The main result of this paper will help us to compute the new invariant.

\section{Topological Background}

In this section we recall some basic definitions and statements from \cite{MoTe1}:\\
Let $D$ be a closed disk on $\R ^2$, $K \subset D$ finite set, $u \in \partial D$. Any diffeomorphism of $D$ which fixes $K$ and is the identity on $\partial D$ acts naturally on $\Pi _1 = \Pi _1(D-K,u)$. We say that two such diffeomorphisms of $D$ (which fix $K$ and equal identity on $\partial D$) are equivalent if they define the same automorphism on $\Pi _1(D-K,u)$. This equivalence relation is compatible with composition of diffeomorphism and thus the equivalence classes form a group.\\


\bde \label{BraidGroupDef}
Braid Group $B_n [D,K]$.
\ede
\noindent Let $D,K$ be as above, $n = \# K$, and let $\mathcal{B}$ be the group of all diffeomorphisms $\beta$ of $D$ such that $\beta (K)=K$, $\beta |_{\partial D} = Id_{\partial D}$. For $\beta _1, \beta _2 \in \mathcal{B}$ we say that $\beta _1$ is equivalent to $\beta _2$ if $\beta _1$ and $\beta _2$ define the same automorphism of $\Pi _1(D-K,u)$. The quotient of $\mathcal{B}$ by this equivalence relation is called the braid group $B_n [D,K]$.\\
Equivalently, if we take the canonical homomorphism $\psi :{\mathcal{B}} \rightarrow {Aut}(\Pi _1(D-K,u))$, then $B_n [D,K] = Im(\psi )$. The elements of $B_n [D,K]$ are called braids.\\ 

\blem
If $K' \subset  D'$ in another pair as above with ${\# K'} = {\# K} = n$, then $B_n [D',K']$ is isomorphic to $B_n [D,K]$.\\
\elem

\noindent This gives rise to the definition of $B_n$:

\bde
$B_n = B_n[D,K]$ for some $D,K$ with $ \# K = n$.
\ede

\bde
$l(q)$.
\ede
\noindent Let $c$ be a small loop equal to the oriented boundary of a small environment $V$ of $w_0$ chosen so that $q' = q - (V \bigcap q)$ is a simple path. Then $l(q) = q' \bigcup c \bigcup {q'}^{-1}$. We also use the notation $l(q)$ for the element of $\Pi _1(D-K,u)$ corresponding to $l(q)$.


\bde
A Bush.
\ede
\noindent Consider in $D$ an ordered set of simple paths $(T_1 ,...,T_n)$ connecting $w_i \in K$ with $u$ such that:\\
\indent (1) $\forall i = 1,...n \quad T_i \bigcap  w_j = \phi $ if $ i \neq j$.\\
\indent (2) $ \bigcap _{i=1}^n T_i = u$.\\ 
\indent (3) For a small circle $c(u)$ around $u$, each ${u_i}' = T_i \bigcap c(u)$ is a single point and the order in $({u'}_1,...{u'}_n)$ is consistent with the positive (``counterclockwise'') orientation of $c(u)$. We say that two such sets $(T_1,...,T_n)$ and $({T'}_1,...,{T'}_n)$ are equivalent if  $\forall i = 1,...,n \quad l(T_i) = l({T'}_i)$ in $\Pi _1(D-K,u)$.\\
An equivalent class of such sets is called a bush in $(D-K,u)$.
The bush represented by $(T_1,...,T_n)$ is denoted by $<T_1,...,T_n>$.\\
It is well known that $\fg$ is a free group.

\bde
g-base
\ede
\noindent A g-base of $\fg$ is an ordered free base of $\fg$ which has the form $(l(T_1),...,l(T_n))$ where $<T_1,...,T_n>$ is a bush in $(D-K,u)$.

\bde \label{HalfTwistDef}
$H(\sigma)$, half-twist defined by $\sigma$.
\ede
\noindent Let $D$ and $K$ be defined as above. Let $a, b \in K$, $K_{a,b} = K \backslash \{ a,b \}$ and $\sigma$ be a simple path (without a self intersection) in $D \backslash \partial D$ connecting $a$ with $b$ such that $\sigma \bigcap K = \{ a,b \}$. Choose a small regular neighborhood  $U$ of $\sigma$ such that $K_{a,b} \bigcap U = \phi$, and an orientation preserving diffeomorphism $\psi : {\R}^2 \rightarrow {\C}^1$ such that
$\psi (\sigma) = [-1,1] = \{ z \in \C ^1 | Re(z) \in [-1,1], Im(z) = 0 \}$ and $\psi (U) = \{ z \in \C ^1||z|<2 \}$. Let $\alpha (r),r \geq 0$, be a real smooth monotone function such that $\alpha (r) = 1$ for $r \in [0,3/2]$ and $\alpha (r) = 0$ for $r \geq 2$. Define a diffeomorphism $h:\C ^1 \rightarrow \C ^1$ as follows: for $z \in \C ^1, z = r e^{i\phi} \quad \text{let} \quad h(z)= r e^{i(\phi +\alpha (r) \pi )}$. It is clear that the restriction of $h$ to $\{ z \in \C ^1 | |z| \leq 3/2 \}$ coincides with the $180 ^\circ$ positive rotation, and that the restriction to $\{ z \in \C ^1 | |z| \geq 2 \}$ is the identity map. The diffeomorphism $\psi ^{-1} \circ h \circ \psi$ induces a braid called half-twist and denoted by $H(\sigma)$\\


\bde \label{frame}
Frame of $B_n[D,K]$
\ede
\noindent Let $K = \{ a_1,...,a_n \}$ and $\sigma _1,..., \sigma _{n-1}$ be a system of simple smooth paths in $D - \partial D$ such that $\sigma _i$ connects $a_i \text{ with } a_{i+1}$ and $L = \bigcup \sigma _i$ is a simple smooth path. The ordered system of half-twists $(H_1,...,H_{n-1})$ defined by $ \{ \sigma _i \} _{i=1} ^{n-1}$ is called a frame of $B_n [D,K]$

\bthm\label{relations} 
Let $(H_1,...,H_{n-1})$ be a frame of $B_n[D,K]$ then,  $B_n[D,K]$ is generated by $\{ H_i \} _{i=1} ^{n-1}$ and the following is a complete list of relations:\\
$H_i H_j = H_j H_i \text{ if }|i-j|>1$ (Commutative relation)\\
$H_i H_j H_i = H_j H_i H_j \text { if } |i-j|=1$ (Triple relation).\\
\ethm

\noindent Proof can be found for example in \cite{MoTe1}.\\
This theorem provides us with Artin's algebraic definition of Braid group. \\

\section{Definition of Hurwitz Moves}


\bde
Hurwitz move on $G^m$ ($R_k, R_k ^{-1}$)
\ede
\noindent Let $G$ be a group, $\overrightarrow{t}=(t_1,...,t_m) \in G^m$. We say that $\overrightarrow{s}=(s_1,...,s_m) \in G^m$ is obtained from $\overrightarrow{t}$ by the Hurwitz move $R_k$ (or $\overrightarrow{t}$ is obtained from $\overrightarrow{s}$ by the Hurwitz move $R_k ^{-1}$) if 
$$ s_i = t_i \quad \text{ for } i \not= k, k+1,$$   
$$ s_k = t_k t_{k+1} t_k ^{-1}, \quad s_{k+1} = t_k.$$
 

\bde
Hurwitz move on factorization 
\ede
\noindent Let $G$ be a group and $t \in G$. Let $t = t_1 \cdots t_m = s_1 \cdots s_m$ be two factorized expressions of $t$. We say that $s_1 \cdots s_m$ is obtained from $t_1 \cdots t_m$ by the Hurwitz move $R_k$ if $(s_1,...,s_m)$ is obtained from $(t_1,...,t_m)$ by the Hurwitz move $R_k$.\\

\bde
Hurwitz equivalence of factorization
\ede
\noindent The factorizations $s_1 \cdots s_m$, $t_1 \cdots t_m$\label{formin2} are Hurwitz equivalent if they are obtained from each other by a finite sequence of Hurwitz moves. The notation is  $t_1 \cdots t_m \overset{HE}{\backsim} s_1 \cdots s_m$.\label{formin1} \\


\bde
Word in $B_n$
\ede
\noindent A word in $B_n$ is a representation of braid as a sequence of the frame elements and their inverses.\\\\


\section{Hurwitz Equivalence of Factorizations with Frame Elements}

\noindent We already proved Proposition \ref{sequence} in \cite{Tz}. We provide here the proof once again since the logic use is basic for understanding the sequel.\\

\noindent We will need a certain result of Garside:
\bcl\label{PositiveEqual}
{\bf (Garside):} Every two positive words (all generators with positive powers) which are equal are transformable into each other through a finite sequence of positive words, such that each word of the sequence is positive and obtained from the preceding one by a direct application of the commutative relation or the triple relation.
\ecl
\noindent Proof: \cite{Garside}.\\

\noindent The following is obvious, nevertheless we give a short proof for clarification:
\bcl \label{commutes}
$G$ is a group $g_1 , g_2 \in G$:\\
1. If $g_1 g_2 = g_2 g_1$ then $g_1  \cdot g_2 \overset{HE}{\backsim} g_2 \cdot g_1$.\\
2. If $g_1 g_2 g_1 = g_2 g_1 g_2$ then $g_1  \cdot g_2 \cdot g_1 \overset{HE}{\backsim} g_2 \cdot g_1 \cdot g_2$.
\ecl
\noindent Proof:\\ 
1. $g_1 \cdot g_2 \overset{R_1}{\rightarrow} g_1 g_2 g_1 ^{-1} \cdot g_1 = g_2 \cdot g_1$.\\
2. $g_1 \cdot g_2 \cdot g_1 \overset{R_2}{\rightarrow} g_1 \cdot g_2 g_1 g_2 ^{-1} \cdot g_2 \overset{R_1}{\rightarrow} g_1 g_2 g_1 g_2 ^{-1} g_1 ^{-1} \cdot g_1 \cdot g_2$ but $g_1 g_2 g_1 = g_2 g_1 g_2$, so $g_1 g_2 g_1 g_2 ^{-1} g_1 ^{-1} = g_2$ and we get $g_1  \cdot g_2 \cdot g_1 \overset{HE}{\backsim} g_2 \cdot g_1 \cdot g_2$\\

\noindent From Garside and the above, we get the following proposition:

\bpr \label{sequence}
Let $H_1,...,H_{n-1}$ be a set of generators of $B_n$ and $H_{i_1} \dots H_{i_p} = H_{j_1} \dots H_{j_p}$ two positive words (with $i_k, j_k \in \{ 1,...n-1 \} $) then  $H_{i_1} \cdots H_{i_p} \overset{HE}{\backsim} H_{j_1} \cdots H_{j_p}$.
\epr

\noindent Proof:\\
Applying \ref{PositiveEqual} on $H_{i_1} \dots H_{i_p} = H_{j_1} \dots H_{j_p}$, we get a finite sequence of positive words $\{ W_r \} _{r=0} ^{q}$ s.t. $W_0 = H_{i_1} \dots H_{i_p}$, $W_q = H_{j_1} \dots H_{j_p}$ and $W_{r+1}$ is obtained from $W_r$ by a single application of the commutative relation or the triple relation.\\
An application of the triple relation (as in \ref{commutes}) is equal to an application of 2 Hurwitz moves $R_{t+1}, R_t$. An application of the commutative relation is equal to an application of one Hurwitz move (on the commuting elements).
Thus, $ W_0 \overset{HE}{\backsim} W_q$

\bde
$\Delta _n ^2  \in B_n[D,K]$
\ede
$\Delta _n ^{2} = (H_1 \dots H_{n-1})^n$ where $ \{ H_i \} _{i=1} ^{n-1}$ is a frame.\\

\noindent We apply \ref{sequence} on $\Delta _n ^2$: 

\bco \label{sameFrameFactorization}
All $\Delta _n ^2$ factorizations $H_{i_1} \cdots H_{i_{n(n-1)}} \quad i_k \in \{ 1,...,n-1\} $ are Hurwitz equivalent.
\eco

\bnote
We use the notation $a[b]$ for conjugating $b$ to $a$.
\enote

\begin{figure}[htp]
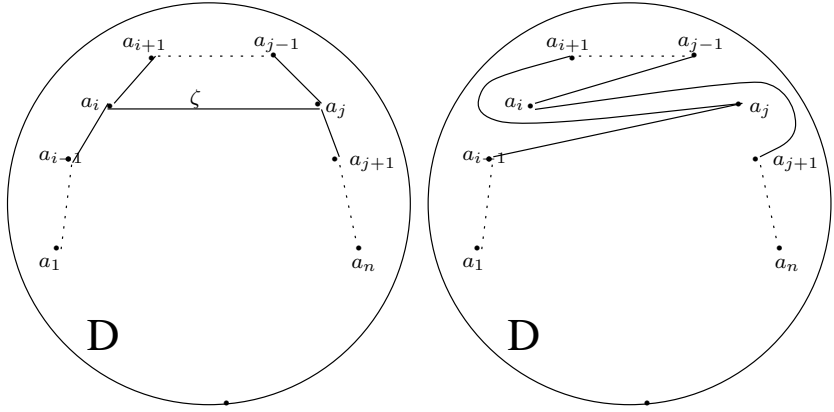

\begin{center}
\epsfxsize=8cm
\epsfysize=4cm
\input conjFrame.pstex_t
\caption[]{
\label{conjFrame}
\sf
$(H_{\zeta}(\sigma _1),...,H_{\zeta}(\sigma _{n-1}))$.}
\end{center}
\end{figure}

\bpr \label{conjFrame1}
Let $\{ H({\sigma _i}) \} _{i=1} ^{n-1}$ be the frame which generates $B_n [D,K]$ and $H_{\zeta}$ a half twist. Then $\{ H({\sigma _i})[H_{\zeta}] \} _{i=1} ^{n-1}$ is also a frame.
\epr
\noindent Proof:\\
By III.1.0 of \cite{MoTe5}, $H(\sigma _i)[H_{\zeta}] = H(H_{\zeta}(\sigma _1))$, thus conjugating $H_{\zeta}$ to the elements in the frame is the same as operating $H_{\zeta}$ on each of the elements in $(\sigma _1,...,\sigma _{n-1})$, and we get $(H_{\zeta}(\sigma _1),...,H_{\zeta}(\sigma _{n-1}))$ which also defines a frame since:\\
\begin{emphit} 
\item $\bigcup \sigma _i$ is a simple smooth path. $H_{\zeta}$ is a diffeomorphism, so $H_{\zeta}(\bigcup \sigma _i) = \bigcup H_{\zeta}(\sigma _i)$ is also a simple smooth path.
\item if $\zeta$ connects $a_i$ with $a_j$, $H_{\zeta}(\sigma _r)$ connects $a_r$ with $a_{r+1}$ when $r \neq i-1,i,j-1,j$. $H_{\zeta}(\sigma _{i-1})$ connects $a_{i-1}$ with $a_j$, $H_{\zeta}(\sigma _{i})$ connects $a_j$ with $a_{i+1}$, $H_{\zeta}(\sigma _{j-1})$ connects $a_{j-1}$ with $a_{i}$, $H_{\zeta}(\sigma _{j})$ connects $a_i$ with $a_{j+1}$. This means that $a_{i-1}$, $a_{i+1}$ are connected through $a_j$, and $a_{j-1}$, $a_{j+1}$ are connected through $a_i$ as shown in Figure \ref{conjFrame}.
\end{emphit}

\bpr \label{moves}
Let $(X_1,...,X_{n-1})$ be the frame which generates $B_n[D,K]$ and $\large{\Pi}$ is the factorization $X_1 \cdots X_{n-1}$, then:\\
1. $X_i \cdot \large{\Pi} \overset{HE}{\backsim} \large{\Pi} \cdot X_{i-1}$ for $1 < i \leq n-1$.\\
2. $X_1 \cdot \large{\Pi} \cdot \large{\Pi} \overset{HE}{\backsim}  \large{\Pi} \cdot \large{\Pi} \cdot X_{n-1}$.
\epr
\noindent Proof:\\
It is enough to show that the words are equal and by \ref{sequence} the factorizations are Hurwitz equivalent.\\
1. $X_i \large{\Pi} = X_i(X_1 ...X_{i-1} X_i ...X_{n-1}) = X_1 ... X_i X_{i-1} X_i...X_{n-1} = \\ 
\quad X_1...X_{i-1} X_i X_{i-1} X_{i+1}...X_{n-1} = (X_1...X_{i-1} X_i X_{i+1}...X_{n-1})X_{i-1} = \large{\Pi} X_{i-1}$.\\
2. $X_1 \large{\Pi} \large{\Pi} = X_1 \large{\Pi} X_1...X_{n-2} X_{n-1}$ and from 1. We get that $ X_1 \large{\Pi} X_1...X_{n-2} X_{n-1} =  X_1 X_2 ... X_{n-1} \large{\Pi} X_{n-1}$.

\bpr \label{oneConj}
Let $ \{ X_i \} _{i=1} ^{n-1}$ be the frame which generates $B_n [D,K]$. For each $ 1 \leq j \leq n-1 $ the $\Delta ^2 _n$ factorizations $ (X_1 \cdots X_{n-1})^n$ and\\ $(X_1 [X_j] \cdots X_{n-1} [X_j])^n$ are Hurwitz equivalent.
\epr
\noindent Proof:\\
We first consider the two extreme cases where $j=n-1$ and $j=1$.\\
In the case where, $j = n-1$, by performing the series of Hurwitz moves $\{ R ^{-1} _{m(n-1)-1} \} _{m=1} ^n$ on $(X_1 \cdots X_{n-1})^n$ (each of the Hurwitz moves is acting on one of the $n$ factors in the product) we get that
$$(X_1 \cdots X_{n-3} \cdot X_{n-2} \cdot X_{n-1})^n \overset{HE}{\backsim} (X_1 \cdots X_{n-3} \cdot X_{n-1} \cdot X_{n-2} [X_{n-1}])^n.$$
Since $X_i [X_{n-1}] = X_i$ for $i \neq {n-2}$, all elements in the factorization are of the form $X_i [X_{n-1}]$ and by \ref{conjFrame1} they are all elements of the same frame. $(X_1 \cdots X_{n-1} \cdot X_{n-2} [X_{n-1}])^n$ and $(X_1 \cdots X_{n-2} [X_{n-1}] \cdot X_{n-1})^n$ are both $\Delta ^2 _n$ factorizations.\\
Having the same frame elements and from \ref{sameFrameFactorization} they are Hurwitz equivalent, and we get that:
$$(X_1 [X_{n-1}] \cdots X_{n-2} [X_{n-1}] \cdot X_{n-1} [X_{n-1}])^n \overset{HE}{\backsim}  (X_1 \cdots X_{n-1})^n.$$
The case where $j=1$ is similar by performing the series of Hurwitz moves $ \{ R ^{-1} _{m(n-1)+1} \} _{m=0} ^{n-1}$.
In the case where $1 < j < n-1$, we perform $\{ R _{m(n-1)+j} ^{-1} \} _{j=0} ^{n-1}$ and get:\\
$$(X_1,\cdots X_{j-1} \cdot X_j \cdot X_{j+1} \cdots X_{n-1})^n \overset{HE}{\backsim} (X_1,\cdots X_{j-2} \cdot X_j \cdot  X_{j-1} [X_j] \cdot X_{j+1} \cdots X_{n-1})^n .$$
Since $X_j$ commutes with $X_1,...,X_{j-2}$, then by \ref{commutes},\\\\
$(X_1,\cdots X_{j-2} \cdot X_j \cdot  X_{j-1} [X_j] \cdot X_{j+1} \cdots X_{n-1})^n \overset{HE}{\backsim} (X_j \cdot X_1,\cdots X_{j-2} \cdot  X_{j-1} [X_j] \cdot X_{j+1} \cdots X_{n-1})^n .$\\\\
(In fact we performed $\{ \{ R_{m(n-1)+k} \} _{k=j-2} ^1 \} _{m=0} ^{n-1}$).

\noindent Again, $X_j$ commutes with $X_{j+2},...,X_{n-1}$ so by \ref{commutes}:\\\\
$(X_j \cdot X_1 \cdots X_{j-2} \cdot  X_{j-1} [X_j] \cdot X_{j+1} \cdots X_{n-1})^n \overset{HE}{\backsim}\\
(X_j \cdot X_1 \cdots X_{j-2} \cdot  X_{j-1} [X_j] \cdot X_{j+1} \cdot X_j \cdot X_{j+2} \cdots X_{n-1}) \cdot  (X_1 \cdots X_{j-2} \cdot  X_{j-1} [X_j] \cdot X_{j+1} \cdot X_j \cdot X_{j+2} \cdots X_{n-1}) ^{n-2} (X_1 \cdots X_{j-2} \cdot  X_{j-1} [X_j] \cdot X_{j+1} \cdot X_{j+2} \cdots X_{n-1})$\\\\
which is Hurwitz equivalent to\\\\
$(X_j \cdot X_1 \cdots X_{j-2} \cdot  X_{j-1} [X_j] \cdot X_j \cdot X_{j+1} [X_j] \cdot X_{j+2} \cdots X_{n-1}) \cdot  (X_1 \cdots X_{j-2} \cdot  X_{j-1} [X_j] \cdot X_j \cdot X_{j+1} [X_j] \cdot X_{j+2} \cdots X_{n-1}) ^{n-2}  (X_1 \cdots X_{j-2} \cdot  X_{j-1} [X_j] \cdot X_{j+1} \cdot X_{j+2} \cdots X_{n-1})$\\\\
(In fact we performed $\{ R ^{-1} _{j+m(n-1)+1} \}_{m=0} ^{n-2}$)\\\\
Let $ (A_1,...,A_{n-1}) $ be the conjugated frame s.t. $A_i = X_i [X_j]$ and $\large{\Pi} = A_1 \cdots A_{n-1}$.\\
The previous factorization can be written as $$A_j \cdot \large{\Pi} ^{n-1} \cdot A_1 \cdots A_{j-1} \cdot X_{j+1} \cdot A_{j+2} \cdots A_{n-1}.$$
by \ref{moves} (1) this factorization is Hurwitz equivalent to: 
$$\large{\Pi} ^{j-1} \cdot A_1 \cdot \large{\Pi} ^{n-j} \cdot A_1 \cdots A_{j-1} \cdot X_{j+1} \cdot A_{j+2} \cdots A_{n-1}.$$
and by \ref{moves} (2) we get that 
$$\large{\Pi} ^{j-1} \cdot A_1 \cdot \large{\Pi} ^{n-j} \cdot A_1 \cdots A_{j-1} \cdot X_{j+1} \cdot A_{j+2} \cdots A_{n-1}$$
is Hurwitz equivalent to
$$\large{\Pi} ^{j+1} \cdot A_{n-1} \cdot \large{\Pi} ^{n-j-2} \cdot A_1 \cdots A_{j-1} \cdot X_{j+1} \cdot A_{j+2} \cdots A_{n-1}$$
Using \ref{moves} (1) $n-j-2$ times we get that the above is Hurwitz equivalent to $$\large{\Pi} ^{n-1} \cdot A_{j+1} \cdot A_1 \cdots A_{j-1} \cdot X_{j+1} \cdot A_{j+2} \cdots A_{n-1}.$$
Since $A_{j+1}$ commutes with $A_1,...,A_{j-1}$ then by \ref{commutes}:\\\\
$\large{\Pi} ^{n-1} \cdot A_{j+1} \cdot A_1 \cdots A_{j-1} \cdot X_{j+1} \cdot A_{j+2} \cdots A_{n-1} \overset{HE}{\backsim}\\ \large{\Pi} ^{n-1} \cdot A_1 \cdots A_{j-1} \cdot A_{j+1} \cdot X_{j+1} \cdot A_{j+2} \cdots A_{n-1}$.\\\\
By performing $R_{(n-1)^2 +j}$ (on $A_{j+1}$), the factorization is Hurwitz equivalent to:
$$ \large{\Pi} ^{n-1} \cdot A_1 \cdots A_{j-1} \cdot X_{j+1} [(X_{j+1}[X_j])^{-1}] \cdot A_{j+1} \cdot A_{j+2} \cdots A_{n-1}$$

\noindent which is exactly $(A_1 \cdots A_{n-1}) ^n$.

\bpr \label{equivalenConj}
Let $b \in B_n[D,K]$, then the factorizations  $(X_1 \cdots X_{n-1})^n$ and $(X_1 [b] \cdots X_{n-1} [b])^n$ are Hurwitz equivalent.
\epr
\noindent Proof:\\
From \cite{Garside}, any braid $P \in B_n[D,K]$ can be represented as $P = \Delta ^r _n \overline{P}$ where $\overline{P} $ is a positive braid. If $r \geq 0$ then clearly $A[P]$ is a conjugation of a positive braid. Assume $r<0$, take $m$ s.t. $2m > -r$ and since $\Delta ^2 _n$ is the generator of the center of $B_n[D,K]$, then $A[P] = A [\Delta ^r _n \overline{P}] = A [\Delta ^{r+2m} _n \overline{P}]$ which is a conjugation of a positive braid.\\
For the rest of the proof we may assume that $P$ is a positive braid.\\
The rest of the proof will be performed by induction on the length of the positive braid $P$. The case where $P=X_{i_1}$ was already proven.\\
If $P = X_{i_1}...X_{i_d} X_{i_{d+1}}$,  by the induction assumption,
$$(X_1 \cdots X_{n-1})^n \overset{HE}{\backsim} (X_1 [X_{i_d}...X_{i_1}] \cdots X_{n-1} [X_{i_d}...X_{i_1}])^n .$$ The new elements are also a frame, so by applying \ref{oneConj} we get the equivalent factorization:
$$((X_1 [X_{i_d}...X_{i_1}])[X_{i_{d+1}}[X_{i_d}...X_{i_1}]] \cdots (X_{n-1} [X_{i_d}...X_{i_1}])[X_{i_{d+1}}[X_{i_d}...X_{i_1}]] )^n .$$
For each $1 \leq j \leq n-1$, $(X_j [X_{i_d}...X_{i_1}])[X_{i_{d+1}}[X_{i_d}...X_{i_1}]] =  X_j  [X_{i_{d+1}}...X_{i_1}]$ and therefore, $$(X_1 \cdots X_{n-1})^n \overset{HE}{\backsim} (X_1 [X_{i_{d+1}}...X_{i_1}] \cdots X_{n-1} [X_{i_{d+1}}...X_{i_1}])^n .$$
\hfill $\qed $

%
%

\section{The Main Result}

\noindent In the following section we will apply the previous proposition on all frames, by proving that every two frames of $B_n[D,K]$ are conjugated.

\begin{figure}[htp]
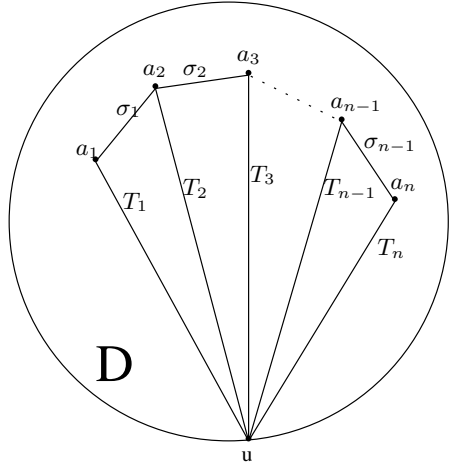

\begin{center}
\epsfxsize=8cm
\epsfysize=4cm
\input gBase.pstex_t
\caption[]{
\label{gBase}
\sf
g-base of a frame.
}
\end{center}
\end{figure}

\bde
g-base of a frame $F$.
\ede
\noindent Let $F=(F_1,...,F_{n-1})$ ($F_i = H(\sigma _i)$) be a frame of $B_n[D,K]$. Let the system of paths $(T_1,...,T_n)$, $T_i$ connecting $u$ with $a_i \in K$ create a bush s.t.\\
$T_i \bigcap \sigma _j =
{	\begin{cases} 
		{\phi,  \text{ if } j \neq i,i-1}\\
		{ a_i, \text{ if } j=i \text{ or } j=i-1}
	\end{cases} } $.\\
Let $\gamma _i = l(T_i )$ then $\Gamma = (\gamma _1 ,...,\gamma _n )$ is the g-base of the frame $F$ as shown in Figure \ref{gBase}.

\begin{figure}[htp]
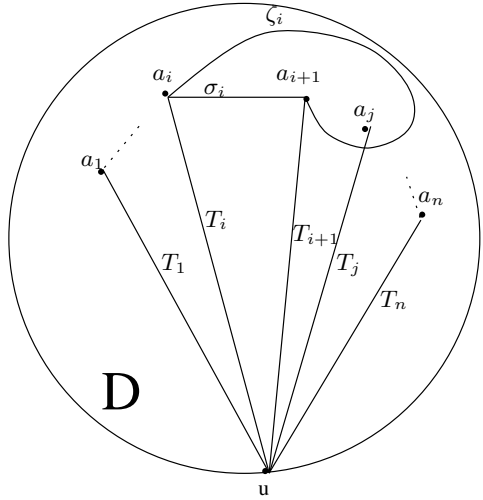

\begin{center}
\epsfxsize=8cm
\epsfysize=4cm
\input twoFrames.pstex_t
\caption[]{
\label{twoFrames}
\sf
Two frames with the same g-base.
}
\end{center}
\end{figure}

\bpr \label{uniqueFrame}
If $\Gamma = (\gamma _1 ,..., \gamma _n)$ is the g-base of two frames, then the frames are equal.
\epr
\noindent Proof:\\
Assuming that $\Gamma$ is the g-base of the two frames $F_1$ and $F_2$ where , $F_1 = (H(\sigma _1),...,H(\sigma _{n-1}))$ and $F_2 = (H(\zeta _1),...,H(\zeta _{n-1}))$ and $F_1 \neq F_2$. Since $F_1 \neq F_2$ then  $\exists i \in \{ 1,...,n-1 \} $ s.t. $H(\sigma _i) \neq H(\zeta _i)$.\\
$\Gamma $ is the g-base of both frames, then $\sigma _i$ and $\zeta _i$ connecting $a_i$ with $a_{i+1}$. $H(\sigma _i) \neq H(\zeta _i)$, therefore, $\zeta _i$ and $\sigma _i$ are ``separated'' by a point $a_j \neq a_i, a_{i+1}$ as shown in Figure \ref{twoFrames}. The two paths start at $a_i$ and end at $a_{i+1}$ therefore creating a cycle with $a_j$ inside. Since $T_j$ connects $a_j$ with $u$, $T_j \bigcap \sigma _i \neq \phi$ or $T_j \bigcap \zeta _i \neq \phi$ and therefore $\Gamma$ is not the g-base of both $F_1$ and $F_2$.

\bpr \label{HTconjFrame}
Let $\Gamma = (\gamma _1 ,..., \gamma _n)$ be the g-base of the frame\\ $F=(H(\sigma _1),...,H(\sigma _{n-1}))$ and let $H_{\zeta}$ be a half twist. Then $(H_{\zeta }(\gamma _1) ,..., H_{\zeta }(\gamma _n) )$ is the g-base of the frame $(H(\sigma _1)[H_{\zeta }],...,H(\sigma _{n-1})[H_{\zeta }])$.
\epr

\noindent Proof:\\
It was already proven in \cite{MoTe1} that any braid $b \in B_n[D,K]$ transforms a g-base to a g-base.
By Proposition \ref{conjFrame1}, $(H_{\zeta}(\sigma _1),...,H_{\zeta}(\sigma _{n-1}))$ defines a frame. 
Since $H_{\zeta}$ is a diffeomorphism, acting on $\sigma _1,...,\sigma _{n-1},\gamma _1,...,\gamma _n$ does not create new intersection points while swapping $a_i$ with $a_j$ for both frame and g-base elements. Therefore, $H_{\zeta}(\Gamma)$ is the g-base of the conjugated frame.

\bco \label{coro1}
Let $\Gamma = (\gamma _1 ,..., \gamma _n)$ be the g-base of $F=(H(\sigma _1),...,H(\sigma _{n-1}))$ and $b \in B_n[D,K]$, then $(b(\gamma _1) ,..., b(\gamma _n) )$ is the g-base of the frame \\ $(H(\sigma _1)[b],...,H(\sigma _{n-1})[b])$.
\eco
\noindent Proof:\\
Induction on the length of $P$ using Proposition \ref{HTconjFrame}.\\

\bpr \label{BconjFrame}
Let $(H(\sigma _1),...,H(\sigma _{n-1}))$ be a frame, then for every frame $(H(\sigma _1 '),...,H(\sigma _{n-1} '))$ there exist a braid $b \in B_n[D,K]$ s.t. $(H(\sigma _1)[b],...,H(\sigma _{n-1})[b]) =(H(\sigma _1 '),...,H(\sigma _{n-1} ')) $.
\epr
\noindent Proof:\\
Let $\Gamma$ and $\Gamma '$ be the g-bases of  $(H(\sigma _1),...,H(\sigma _{n-1}))$ and  $(H(\sigma _1 '),...,H(\sigma _{n-1} '))$.\\
In \cite{MoTe1} it was proved that there exists a unique braid $b \in B_n[D,K]$ s.t. $b(\Gamma )= \Gamma ' $. From Corollary \ref{coro1} $(H(\sigma _1)[b],...,H(\sigma _{n-1})[b])$ is a frame with a g-base $\Gamma '$.\\
From Proposition \ref{uniqueFrame}, $\Gamma '$ has only one frame, therefore, $(H(\sigma _1)[b],...,H(\sigma _{n-1})[b]) = (H(\sigma _1 '),...,H(\sigma _{n-1} '))$\\

\noindent We now come to the main result:

\bthm
Let $(H_1,...,H_{n-1})$ and $(F_1,...,F_{n-1})$ be two frames of $B_n [D,K]$, then for every two braids s.t. $\Delta _n ^2 = F_{i_1}...F_{i_{n(n-1)}} = H_{j_1}...H_{j_{n(n-1)}}$ the factorizations $ F_{i_1} \cdots F_{i_{n(n-1)}}$ and $H_{j_1} \cdots H_{j_{n(n-1)}}$ are Hurwitz equivalent.
\ethm
\noindent Proof:\\
From Corollary \ref{sameFrameFactorization}, $ F_{i_1} \cdots F_{i_{n(n-1)}} \overset{HE}{\backsim} (F_1 \cdots F_{n-1})^n$ and $H_{j_1} \cdots H_{j_{n(n-1)}} \overset{HE}{\backsim} (H_1 \cdots H_{n-1})^n$.\\
From Proposition \ref{BconjFrame}, $\exists b \in B_n[D,K]$ s.t. $(F_1,...,F_{n-1}) = (H_1 [b],...H_{n-1} [b])$, so from Proposition \ref{equivalenConj}, $(F_1 \cdots F_{n-1})^n \overset{HE}{\backsim} (H_1 \cdots H_{n-1})^n$ and therefore, $ F_{i_1} \cdots F_{i_{n(n-1)}}, \quad H_{j_1} \cdots H_{j_{n(n-1)}}$ are Hurwitz equivalent. $\hfill  \qed$

\end{document}

%% file: conjFrame.pstex_t
\begin{picture}(0,0)%
\epsfig{file=conjFrame.pstex}%
\end{picture}%
\setlength{\unitlength}{0.00083300in}%
\begingroup\makeatletter\ifx\SetFigFont\undefined
\def\x#1#2#3#4#5#6#7\relax{\def\x{#1#2#3#4#5#6}}%
\expandafter\x\fmtname xxxxxx\relax \def\y{splain}%
\ifx\x\y   
\gdef\SetFigFont#1#2#3{%
  \ifnum #1<17\tiny\else \ifnum #1<20\small\else
  \ifnum #1<24\normalsize\else \ifnum #1<29\large\else
  \ifnum #1<34\Large\else \ifnum #1<41\LARGE\else
     \huge\fi\fi\fi\fi\fi\fi
  \csname #3\endcsname}%
\else
\gdef\SetFigFont#1#2#3{\begingroup
  \count@#1\relax \ifnum 25<\count@\count@25\fi
  \def\x{\endgroup\@setsize\SetFigFont{#2pt}}%
  \expandafter\x
    \csname \romannumeral\the\count@ pt\expandafter\endcsname
    \csname @\romannumeral\the\count@ pt\endcsname
  \csname #3\endcsname}%
\fi
\fi\endgroup
\begin{picture}(5188,2540)(1418,-3419)
\put(4526,-1543){\makebox(0,0)[lb]{\smash{\SetFigFont{8}{9.6}{rm}$a_i$}}}
\put(6054,-1573){\makebox(0,0)[lb]{\smash{\SetFigFont{8}{9.6}{rm}$a_j$}}}
\put(4265,-2548){\makebox(0,0)[lb]{\smash{\SetFigFont{8}{9.6}{rm}$a_1$}}}
\put(6214,-2553){\makebox(0,0)[lb]{\smash{\SetFigFont{8}{9.6}{rm}$a_n$}}}
\put(4786,-1171){\makebox(0,0)[lb]{\smash{\SetFigFont{8}{9.6}{rm}$a_{i+1}$}}}
\put(5607,-1163){\makebox(0,0)[lb]{\smash{\SetFigFont{8}{9.6}{rm}$a_{j-1}$}}}
\put(4265,-1878){\makebox(0,0)[lb]{\smash{\SetFigFont{8}{9.6}{rm}$a_{i-1}$}}}
\put(6202,-1908){\makebox(0,0)[lb]{\smash{\SetFigFont{8}{9.6}{rm}$a_{j+1}$}}}
\put(1883,-1543){\makebox(0,0)[lb]{\smash{\SetFigFont{8}{9.6}{rm}$a_i$}}}
\put(3411,-1573){\makebox(0,0)[lb]{\smash{\SetFigFont{8}{9.6}{rm}$a_j$}}}
\put(1623,-2548){\makebox(0,0)[lb]{\smash{\SetFigFont{8}{9.6}{rm}$a_1$}}}
\put(3572,-2553){\makebox(0,0)[lb]{\smash{\SetFigFont{8}{9.6}{rm}$a_n$}}}
\put(2144,-1171){\makebox(0,0)[lb]{\smash{\SetFigFont{8}{9.6}{rm}$a_{i+1}$}}}
\put(2964,-1163){\makebox(0,0)[lb]{\smash{\SetFigFont{8}{9.6}{rm}$a_{j-1}$}}}
\put(1623,-1878){\makebox(0,0)[lb]{\smash{\SetFigFont{8}{9.6}{rm}$a_{i-1}$}}}
\put(3560,-1908){\makebox(0,0)[lb]{\smash{\SetFigFont{8}{9.6}{rm}$a_{j+1}$}}}
\put(2564,-1528){\makebox(0,0)[lb]{\smash{\SetFigFont{8}{9.6}{rm}$\zeta$}}}
\end{picture}

%% file: gBase.pstex_t
\begin{picture}(0,0)%
\epsfig{file=gBase.pstex}%
\end{picture}%
\setlength{\unitlength}{0.00083300in}%
\begingroup\makeatletter\ifx\SetFigFont\undefined
\def\x#1#2#3#4#5#6#7\relax{\def\x{#1#2#3#4#5#6}}%
\expandafter\x\fmtname xxxxxx\relax \def\y{splain}%
\ifx\x\y   
\gdef\SetFigFont#1#2#3{%
  \ifnum #1<17\tiny\else \ifnum #1<20\small\else
  \ifnum #1<24\normalsize\else \ifnum #1<29\large\else
  \ifnum #1<34\Large\else \ifnum #1<41\LARGE\else
     \huge\fi\fi\fi\fi\fi\fi
  \csname #3\endcsname}%
\else
\gdef\SetFigFont#1#2#3{\begingroup
  \count@#1\relax \ifnum 25<\count@\count@25\fi
  \def\x{\endgroup\@setsize\SetFigFont{#2pt}}%
  \expandafter\x
    \csname \romannumeral\the\count@ pt\expandafter\endcsname
    \csname @\romannumeral\the\count@ pt\endcsname
  \csname #3\endcsname}%
\fi
\fi\endgroup
\begin{picture}(2774,2902)(1385,-3747)
\put(1811,-1814){\makebox(0,0)[lb]{\smash{\SetFigFont{9}{10.8}{rm}$a_1$}}}
\put(2228,-1312){\makebox(0,0)[lb]{\smash{\SetFigFont{9}{10.8}{rm}$a_2$}}}
\put(2813,-1229){\makebox(0,0)[lb]{\smash{\SetFigFont{9}{10.8}{rm}$a_3$}}}
\put(3398,-1521){\makebox(0,0)[lb]{\smash{\SetFigFont{9}{10.8}{rm}$a_{n-1}$}}}
\put(3774,-2023){\makebox(0,0)[lb]{\smash{\SetFigFont{9}{10.8}{rm}$a_n$}}}
\put(2103,-2148){\makebox(0,0)[lb]{\smash{\SetFigFont{9}{10.8}{rm}$T_1$}}}
\put(2479,-2065){\makebox(0,0)[lb]{\smash{\SetFigFont{9}{10.8}{rm}$T_2$}}}
\put(2897,-1981){\makebox(0,0)[lb]{\smash{\SetFigFont{9}{10.8}{rm}$T_3$}}}
\put(3357,-2065){\makebox(0,0)[lb]{\smash{\SetFigFont{9}{10.8}{rm}$T_{n-1}$}}}
\put(3691,-2441){\makebox(0,0)[lb]{\smash{\SetFigFont{9}{10.8}{rm}$T_n$}}}
\put(2061,-1563){\makebox(0,0)[lb]{\smash{\SetFigFont{9}{10.8}{rm}$\sigma _1$}}}
\put(3607,-1772){\makebox(0,0)[lb]{\smash{\SetFigFont{9}{10.8}{it}$\sigma _{n-1}$}}}
\put(2479,-1312){\makebox(0,0)[lb]{\smash{\SetFigFont{9}{10.8}{rm}$\sigma _2$}}}
\end{picture}

%% file: twoFrames.pstex_t
\begin{picture}(0,0)%
\epsfig{file=twoFrames.pstex}%
\end{picture}%
\setlength{\unitlength}{0.00083300in}%
\begingroup\makeatletter\ifx\SetFigFont\undefined
\def\x#1#2#3#4#5#6#7\relax{\def\x{#1#2#3#4#5#6}}%
\expandafter\x\fmtname xxxxxx\relax \def\y{splain}%
\ifx\x\y   
\gdef\SetFigFont#1#2#3{%
  \ifnum #1<17\tiny\else \ifnum #1<20\small\else
  \ifnum #1<24\normalsize\else \ifnum #1<29\large\else
  \ifnum #1<34\Large\else \ifnum #1<41\LARGE\else
     \huge\fi\fi\fi\fi\fi\fi
  \csname #3\endcsname}%
\else
\gdef\SetFigFont#1#2#3{\begingroup
  \count@#1\relax \ifnum 25<\count@\count@25\fi
  \def\x{\endgroup\@setsize\SetFigFont{#2pt}}%
  \expandafter\x
    \csname \romannumeral\the\count@ pt\expandafter\endcsname
    \csname @\romannumeral\the\count@ pt\endcsname
  \csname #3\endcsname}%
\fi
\fi\endgroup
\begin{picture}(2974,3108)(1418,-3987)
\put(1874,-1915){\makebox(0,0)[lb]{\smash{\SetFigFont{9}{10.8}{rm}$a_1$}}}
\put(2321,-1377){\makebox(0,0)[lb]{\smash{\SetFigFont{9}{10.8}{rm}$a_i$}}}
\put(3574,-1602){\makebox(0,0)[lb]{\smash{\SetFigFont{9}{10.8}{rm}$a_j$}}}
\put(3979,-2139){\makebox(0,0)[lb]{\smash{\SetFigFont{9}{10.8}{rm}$a_n$}}}
\put(3093,-1377){\makebox(0,0)[lb]{\smash{\SetFigFont{9}{10.8}{rm}$a_{i+1}$}}}
\put(2642,-1465){\makebox(0,0)[lb]{\smash{\SetFigFont{9}{10.8}{rm}$\sigma _i$}}}
\put(2379,-2566){\makebox(0,0)[lb]{\smash{\SetFigFont{9}{10.8}{rm}$T_1$}}}
\put(2651,-2293){\makebox(0,0)[lb]{\smash{\SetFigFont{9}{10.8}{rm}$T_i$}}}
\put(3198,-2384){\makebox(0,0)[lb]{\smash{\SetFigFont{9}{10.8}{rm}$T_{i+1}$}}}
\put(3470,-2566){\makebox(0,0)[lb]{\smash{\SetFigFont{9}{10.8}{rm}$T_j$}}}
\put(3742,-2792){\makebox(0,0)[lb]{\smash{\SetFigFont{9}{10.8}{rm}$T_n$}}}
\put(3027,-1008){\makebox(0,0)[lb]{\smash{\SetFigFont{9}{10.8}{rm}$\zeta _i$}}}
\end{picture}

%% file: b4_1.bbl
\newpage

\begin{thebibliography}{99}
\bibitem{Garside} {F. A. Garside, {\it The Braid group and other groups} Quart. J. Math Oxford (2) 78 (1969), 235-254.}
\bibitem{MoTe5} {B. Moishezon and M. Teicher, {\it Braid group technique in complex geometry V: The fundamental group of a complement of a branch curve of a veronese generic projection}, Communication in Analysis and Geometry 4, No. 11, (1996), 1-120}
\bibitem{MoTe1} {B. Moishezon and M. Teicher, {\it Braid group techniques in complex geometry I, Line arrangements in $\C {\Bbb P} ^2$}, Contemporary Math. 78 (1988), 425-555.}
\bibitem{KuTe} {Vik. S. Kulikov and M. Teicher, {\it  Braid monodromy factorization and diffeomorphism types}, Izvestia, Journal of the Russian Academy of Science, Tom 64, 2, 2000, 89-120.}
\bibitem{Te} {M. Teicher, {\it Braid monodromy type of 4 manifolds}, Preprint.}
\bibitem{Tz} {M. Teicher and T. Ben-Itzhak, {\it Hurwitz equivalence in braid group $B_3$}, Preprint.}

\end{thebibliography}
